\documentclass[leqno,a4paper,12pt]{article}
\usepackage{amsfonts,amssymb}

\def\Z{\mathbb{Z}}
\def\Q{\mathbb{Q}}
\def\R{\mathbb{R}}
\def\C{\mathbb{C}}
\def\F{\mathbb{F}}
\def\G{\mathbb{G}}
\def\P{\mathbb{P}}

\def\A{\mathbb{A}}
\def\O{{\cal O}}

\def\til#1{\widetilde{#1}}
\def\ovl#1{\overline{#1}}
\def\pf{{\indent\textsc{Proof.}\ }}
\def\qed{\hfill$\square$}
\def\GL{\mathrm{GL}}
\def\PGL{\mathrm{PGL}}
\def\M{\mathrm{M}}
\def\Spec{\mathop{\mathrm{Spec}}\nolimits}
\def\Spf{\mathop{\mathrm{Spf}}\nolimits}
\def\tr{\mathop{\mathrm{tr}}\nolimits}
\def\det{\mathop{\mathrm{det}}\nolimits}
\def\hasseinv{\mathop{\mathrm{inv}}\nolimits}
\def\opp#1{{{#1}}^{\mathrm{op}}}
\def\inv#1{#1^{\times}}
\def\ad#1{{{#1}}_{\mathrm{ad}}}
\begin{document}
\newcounter{state}[section]
\setcounter{state}{0}
\setcounter{section}{0}
\renewcommand{\thestate}{\arabic{section}.\arabic{state}}
\renewcommand{\thesection}{\arabic{section}}
\newcommand{\statenumber}{\refstepcounter{state}
\indent{\bf\thestate}}
\newcommand{\sectionnumber}{\refstepcounter{section}
\indent{\bf \thesection.}}
\begin{center}
\begin{large}
{\sf Arithmetic structure of Mumford's fake projective plane}
\footnote{{\em $1991$ Mathematics Subject Classification\/}.
Primary 14G35; Secondary 14J29, 11G15, 11Q25.}
\end{large}

\vspace{4ex}
{\sc Fumiharu Kato}

\vspace{4ex}
\begin{minipage}{10.5cm}
\setlength{\baselineskip}{.85\baselineskip}
\begin{small}
We prove that Mumford's fake projective plane is a connected Shimura
variety associated to a certain unitary group.
We also describe explicitly the unitary group and the necessary data for 
defining this Shimura variety, and give a field of definition.
\end{small}
\end{minipage}
\end{center}

\vspace{1ex}
\sectionnumber\label{introduction}\ {\bf Introduction}

\vspace{1ex}
In this paper we construct explicitly a Shimura surface of PEL-type
whose connected components are, as complex surfaces, isomorphic to
Mumford's fake projective plane.
A fake projective plane is a non-singular projective surface of general
type having the same topological Betti numbers as $\C\P^2$.
The first example of fake projective planes was discovered by Mumford in
\cite{Mu79}, where he constructed it by using the theory of $p$-adic
uniformization. 
Later in part of \cite{IK98} Masanori Ishida and the author observed 
other possible fake projective planes according to Mumford's idea, and
proved there exist (at least) two more different fake projective planes.
As far as the author knows, these three are the only known fake
projective planes.

A fake projective plane has complex unit-ball as the topological universal
covering, and in this sense it is among the most interesting classes of
algebraic surfaces. 
However, very little is known about its algebraic or analytic structure
(recently Barlow \cite{Ba99} investigated zero-cycles on Mumford's fake
projective plane in connection with Bloch conjecture).
For instance, topological or even some algebro-geometric properties
of them (e.g.\ fundamental groups, field of definition, etc.) are
quite obscure; even the arithmeticity of the fundamental groups seems
unknown.
The main reason for this is that the known fake projective planes are, a
priori, constructed as surfaces over the field of $2$-adic numbers
$\Q_2$, not over $\C$, and hence are not quite concrete as complex surfaces.
(A fake projective plane over an arbitrary field may be defined to be a
smooth surface enjoying: $\mathrm{c}^2_1=3\mathrm{c}_2=9$, 
$\mathrm{p}_{\mathrm{g}}=\mathrm{q}=0$ and the canonical class is
ample; if one considers this condition in the complex analytic category,
this condition is equivalent to that the surface in question is a fake
projective plane in the origial sense.)

The theory of $p$-adic uniformization, on the other hand, involves
surprizingly rich algebro-geometric and arithmetic aspects concerning 
with Shimura varieties. 
This was first observed by \v{C}erednik \cite{Ce76} and Drinfeld \cite{Dr76} 
for Shimura curves.
The generalizations to higher dimensions have been developed over the past
two decades; recently, Rapoport-Zink \cite{RZ96}, Boutot-Zink
\cite{BZ95} and Varshavsky \cite{Va98a}\cite{Va98b} gave powerfully
extended generalities. 
These theories give a nice symmetric picture involving $p$-adic and
complex uniformizations of certain Shimura varieties of
PEL-type; viz.\ such a Shimura variety has both $p$-adic (for very
special $p$) and complex analytic uniformizations, and they can be viewed
symmetrically by means of ``changing-invariants'' trick.
Shimura varieties of this kind are therefore much expected to possess 
rich structures involving the interplay between $p$-adic and complex
analysis.
However, as far as the author knows, very few such examples are known in
dimension two or high, in contract with numerous explicit examples of Shimura 
curves of this type (see, e.g.\ \cite{Ku79}\cite{Ku94}\cite{An98}).

\vspace{1ex}
{\it Sketch of Results.}\ 
In the next section we will construct explicitly a Shimura variety 
${\cal S}h_C$ associated to a certain central division algebra over
$K=\Q(\sqrt{-7})$. 
The Shimura variety ${\cal S}h_C$ has the canonical model $\mathrm{Sh}_C$
over $K$.
We will apply the theory by Rapoport-Zink \cite{RZ96}, Boutot-Zink
\cite{BZ95} and Varshavsky \cite{Va98a}\cite{Va98b} to this Shimura
variety to show the following results (Theorem \ref{thm-main} and
Corollary \ref{cor-main}): 

\vspace{1ex}
{\bf Theorem.}\ {\sl The canonical model
$\mathrm{Sh}_C$ of the Shimura variety ${\cal S}h_C$ has exactly three
geometrically connected components defined over $\Q(\zeta_7)$ which are
permuted by the action of $\mathrm{Gal}(\Q(\zeta_7)/\Q(\sqrt{-7}))$.
Moreover, for any connected component $S$ of $\mathrm{Sh}_C$ there
exists an isomorphism 
$$
X_{\mathrm{Mum}}\times_{\Spec\Q_2}\Spec\Q_2(\zeta_7)
\stackrel{\sim}{\longrightarrow} 
S\times_{\Spec\Q(\zeta_7)}\Spec\Q_2(\zeta_7)
$$
over $\Q_2(\zeta_7)$, where $X_{\mathrm{Mum}}$ is the Mumford's fake
projective plane.}

\vspace{1ex}
{\bf Corollary.}\ {\sl 
Each connected component of the complex surface ${\cal S}h_C$ is a fake
projective plane, 
which has a field of definition $\Q(\zeta_7)$.
Moreover, it is an arithmetic quotient of complex
unit-ball.}

\vspace{1ex}
The structure of this paper is as follows: In the next section we
construct the Shimura variety ${\cal S}h$ by presenting the necessary
data explicitly. 
In \S\ref{section-construction} we review briefly Mumford's construction
of the fake projective plane. 
In the final section \S\ref{section-relation} we investigate the
relation between the fake projective plane and the Shimura variety
${\cal S}h$, and prove the main theorem.

\vspace{1ex}
{\it Notation and conventions.}\ 
The notation $A\otimes B$ (absence of the base ring) only occurs when 
$A$ and $B$ are $\Q$-algebra, and the tensor is taken over $\Q$.
For a field extension $F/E$, we denote by $\mathrm{Res}_{F/E}$ the Weil
restriction. For a ring $A$ we denote by $\opp{A}$ the opposite ring. 
By an involution of a ring $A$ we always mean a homomorphism
$\ast\colon A\rightarrow\opp{A}$ such that $\ast\circ\ast=\mathrm{id}_A$.

\vspace{2ex}
\sectionnumber\label{section-construction}\ 
{\bf Construction of Shimura variety}

\vspace{1ex}
\statenumber\label{para-algebra}{\bf .}\ 
In this section we construct a Shimura variety which will later be
related to Mumford's fake projective plane. Our Shimura variety ${\cal
S}h$ is of PEL-type related to a certain unitary group. In order to
define it we introduce 

\vspace{1ex}
$\bullet$ a central division algebra $D$ over $K$,

$\bullet$ a positive involution $\ast$ of $D$ of second kind,

$\bullet$ and a non-degenerate anti-hermitian form $\psi$.

\vspace{1ex}\noindent
These data yield the Shimura data (cf.\ \cite{De71}) of our Shimura
variety.

Let $\zeta=\zeta_7$ be a primitive $7$th-root of unity, and
set $L=\Q(\zeta)$; $L$ is a cyclic extension of $\Q$ of degree $6$
having the intermediate quadratic extension $K=\Q(\lambda)\ (\cong
\Q(\sqrt{-7}))$, where $\lambda=\zeta+\zeta^2+\zeta^4$.
The Galois group of $L/K$ is generated by the Frobenius map $\sigma\colon
\zeta\mapsto\zeta^2$, and that of $K/\Q$ is generated by the complex 
conjugation $z\mapsto\ovl{z}$.
Note that the prime $2$ decomposes on $K$ such that 
$2=\lambda\ovl{\lambda}$ gives the prime factorization.
The prime $7$ ramifies on $K$.
We fix an infinite place (CM-type) $\varepsilon\colon K\hookrightarrow\C$ by
$\lambda\mapsto\frac{-1+\sqrt{-7}}{2}$.
Note that the class number of $K$ is $1$.
We set 
$$
\mu=\lambda/\ovl{\lambda}.
$$
Now the central division algebra $D$ over $K$ (of dimension $9$) is
defined by 
$$
D=\bigoplus^2_{i=0}L\Pi^i;\quad\Pi^3=\mu,\ \Pi z=z^{\sigma}\Pi\ \ 
(z\in L).
$$

{\bf \ref{para-algebra}.1.}\ {\bf Lemma.}\ 
{\sl The $L$-algebra $D\otimes_KL$ is isomorphic to the matrix algebra 
$\M_3(L)$. In particular, we have $D_{\varepsilon}\cong\M_3(\C)$.} 

\vspace{1ex}
This is well-known; for the later purpose we give an isomorphism
$\Phi\colon D\otimes_KL\stackrel{\sim}{\rightarrow}\mathrm{M}_3(L)$
explicitly as follows: 
Set $V=D$, which we consider as a left $D$-module.
Then we have
$$
V\otimes_KL\stackrel{\sim}{\longrightarrow}
\bigoplus^2_{i=0}V\otimes_{L,\sigma^i}L.
$$
Looking at the action of $D$ on the factor $V\otimes_{L,\sigma^0}L$
subject to the base $\{1\otimes 1,\Pi\otimes 1,\Pi^2\otimes 1\}$ we
obtain a matrix presentation of elements in $D$ determined by 
$$
z\mapsto\left[
\begin{array}{ccc}
z\\
&z\rlap{$\,\!^{\sigma^2}$}\\
&&z\rlap{$\,\!^{\sigma}$}
\end{array}
\right]
\quad\textrm{and}\quad
\Pi\mapsto\left[
\begin{array}{ccc}
&&\mu\\
1\\
&1
\end{array}
\right],
\leqno{(\ref{para-algebra}.2)}
$$
which induces the desired isomorphism $\Phi$.
By taking an embedding $\theta\colon L\hookrightarrow\C$ by
$\zeta\mapsto\mathrm{exp}\frac{2i\pi}{7}$ we can extend $\Phi$ to an
isomorphism $D_{\varepsilon}\stackrel{\sim}{\rightarrow}\M_3(\C)$.

\vspace{1ex}
{\bf \ref{para-algebra}.3.}\ {\bf Lemma.}\ 
{\sl We have $\hasseinv_{\lambda}D=1/3$ and
$\hasseinv_{\ovl{\lambda}}D=-1/3$, while for any finite place $\ell$
which does not divide $2$ we have $\hasseinv_{\ell}D=0$.}

\vspace{1ex}
This is clear by the construction.
We define the order $\O_D$ of $D$ by
$$
\O_D=\O_L\oplus\O_L\ovl{\lambda}\Pi\oplus\O_L\ovl{\lambda}\Pi^2.
\leqno{(\ref{para-algebra}.4)}
$$
It is easy to verify that
$\O_{D_{\lambda}}=\O_D\otimes_{\O_K}\O_{K_{\lambda}}$ and
$\O_{D_{\ovl{\lambda}}}=\O_D\otimes_{\O_K}\O_{K_{\ovl{\lambda}}}$ are
maxial orders of $D_{\lambda}=D\otimes_KK_{\lambda}$ and 
$D_{\ovl{\lambda}}=D\otimes_KK_{\ovl{\lambda}}$, respectively.

\vspace{1ex}
\statenumber\label{para-involution}{\bf .}\ Next we define an involution
(positive of second kind) $\ast$ on $D$ by
$$
\Pi^{\ast}=\ovl{\mu}\Pi^2\quad\textrm{and}\quad
z^{\ast}=\ovl{z}\ \ (z\in L).
$$
It is easy to check that the involution $\ast$ is mapped by the
isomorphism $\Phi$ as above to the standard involution $(a_{ij})^{\ast}=
(\ovl{a}_{ji})$. 
Set 
$$
b=(\lambda-\ovl{\lambda})-\ovl{\lambda}\Pi+\ovl{\lambda}\Pi^2.
\leqno(\ref{para-involution}.1)
$$
The element $b$ is anti-symmetric (i.e.\ $b^{\ast}=-b$); we have
$$
\Phi(b)=\left[
\begin{array}{ccc}
\lambda-\ovl{\lambda}&\lambda&\llap{$-$}\lambda\\
\llap{$-$}\ovl{\lambda}&\lambda-\ovl{\lambda}&\lambda\\
\ovl{\lambda}&\llap{$-$}\ovl{\lambda}&\lambda-\ovl{\lambda}
\end{array}
\right].
\leqno(\ref{para-involution}.2)
$$
The element $b$ induces an involution, positive of second kind,
$\bigstar\colon D\rightarrow\opp{D}$ by
$\alpha^{\bigstar}=b\alpha^{\ast}b^{-1}$ and a $\Q$-bilinear form
$\psi\colon V\times V\rightarrow\Q$ on $V$ by 
$$
\psi(x,y)=\tr_{D/\Q}(ybx^{\ast})
$$
for $x,y\in V$.

\vspace{1ex}
\statenumber\label{lem-signature}{\bf .\ Lemma.}\ 

1) {\sl The bilinear form $\psi$ is non-degenerate and
anti-symmetric. Moreover we have $\psi(\alpha
x,y)=\psi(x,\alpha^{\ast}y)$ for $\alpha\in D$ and $x,y\in V$.} 

2) {\sl There exists $P\in\M_3(\C)$ such that $P^{\ast}\Phi(b)P=
\mathrm{diag}(-\sqrt{-1},-\sqrt{-1},\sqrt{-1})$.}

\vspace{1ex}
\pf
1) Non-degeneracy is clear since $b$ is invertible. The
anti-symmetricity follows from the equality $b=-b^{\ast}$. The other
assertion is clear. 
The characteristic polynomial of $\Phi(b)$ is
$t^3-3\sqrt{-7}t^2-15t-\sqrt{-7}$. By this 2) follows immediately.
\qed

\vspace{1ex}
\statenumber\label{def-group}{\bf .\ Definition.}\ 
We define an algebraic group
$G$ over $\Q$ by
\begin{eqnarray*}
G(R)&=&\{\gamma\in\inv{(D\otimes R)}\ |\ \psi(\gamma x,\gamma y)=
c(\gamma)\psi(x,y),\ c(\gamma)\in\inv{R},\ x,y\in V\otimes R\}\\
&=&\{\gamma\in\inv{(D\otimes R)}\ |\ \gamma\gamma^{\bigstar}\in R\}
\end{eqnarray*}
for any (commutative) $\Q$-algebra $R$.

\vspace{1ex}
{\bf \ref{def-group}.1.}\ {\bf Lemma.}\ {\sl The isomorphism $\Phi\colon 
D\otimes\R\stackrel{\sim}{\rightarrow}\M_3(\C)$ defined by}
(\ref{para-algebra}.2) {\sl induces an isomorphism
$\Phi\colon G_{\R}\stackrel{\sim}{\rightarrow}\mathrm{GU}(2,1)$.} 

\vspace{1ex}
This is clear by \ref{lem-signature}.2.

\vspace{1ex}
\statenumber\label{para-shimura}{\bf .}\ Consider the homomorphism of
algebraic groups over $\R$
$$
h\colon\mathrm{Res}_{\C/\R}\G_{\mathrm{m}}\longrightarrow G_{\R}
\leqno{(\ref{para-shimura}.1)}
$$ 
such that $\Phi\circ
h(\sqrt{-1})=\mathrm{diag}(\sqrt{-1},\sqrt{-1},-\sqrt{-1})$ and 
$\Phi\circ h(r)$ is the scalar matrix $r$ for $r\in\inv{\R}$.
By \cite[4.1]{Ko92} the pair $(G_{\R},h)$ satisfies the conditions
(1.5.1), (1.5.2) and (1.5.3) in \cite{De71}, and hence we get the 
associated Shimura variety;
for a sufficiently small open compact subgroup $C$ of
$G(\A_{\mathrm{f}})$ the corresponding Shimura variety ${\cal S}h_C$ is 
defined by
$$
{\cal S}h_C=G(\Q)\setminus X_{\infty}\times 
G(\A_{\mathrm{f}})/C,
\leqno{(\ref{para-shimura}.2)}
$$
where $X_{\infty}$ is the set of all conjugates of $h$ under the action
of $G(\R)$.
Note that $X_{\infty}$ is isomorphic to the open unit-ball in $\C^2$.
The Shimura variety ${\cal S}h_C$ is a finite disjoint union of
quasi-projective manifolds, which are arithmetic quotients of the complex
unit-ball.
It is easy to see that the Shimura field $E$ coincides with
$\varepsilon(K)$. 
We therefore obtain the canonical model $\mathrm{Sh}_C$, which is a
quasi-projective variety over $E=\Q(\sqrt{-7})$.

\vspace{2ex}
\sectionnumber\label{section-mumford}\ 
{\bf Mumford's fake projective plane}

\vspace{1ex}
\statenumber\label{para-recall}{\bf .}\ 
In this section we briefly review the Mumford's construction of a fake
projective plane. 
We mainly focus on the construction of the discrete group
$\Gamma_{\mathrm{Mum}}$ which gives the uniformization of the fake
projective plane. 
For more details the reader should consult Mumford's original paper
\cite{Mu79}. 

We use the notation as in \ref{para-algebra}.
Regarding $L$ as a $3$-dimensional $K$-vector space we define a
hermitian form $h\colon L\times L\rightarrow\Q$ by
$h(x,y)=\tr_{L/K}(x\ovl{y})$ for $x,y\in L$.
Subject to the basis $1,\zeta,\zeta^2$ of $L$ the form $h$ is
represented by the matrix
$$
H=\left[
\begin{array}{lll}
3&\ovl{\lambda}&\ovl{\lambda}\\
\lambda&3&\ovl{\lambda}\\
\lambda&\lambda&3
\end{array}
\right],
\leqno{(\ref{para-recall}.1)}
$$
which is a positive definite hermitian matrix of determinant $7$.
Note that there exists a decomposition over $\Q$
$$
H=WW^{\ast}\quad\textrm{with}\quad
W=\left[
\begin{array}{ccc}
\lambda&1&0\\
0&\lambda&1\\
\mu&0&\lambda
\end{array}
\right].
\leqno{(\ref{para-recall}.2)}
$$
Let $\M_3(K)$ be the set of all $3\times 3$ matrices with entries in
$K$, regarded as an algebra over $\Q$.
We denote by $\ast\colon\M_3(K)\rightarrow\opp{\M_3(K)}$ the standard
involution $(a_{ij})^{\ast}=(\ovl{a}_{ji})$.
We define another involution (positive of second kind)
$\dagger\colon\M_3(K)\rightarrow\opp{\M_3(K)}$ by
$A^{\dagger}=HA^{\ast}H^{-1}$.

\vspace{1ex}
\statenumber\label{def-twistgroup}{\bf .\ Definition.}\ 
We define an algebraic group
$I$ over $\Q$ by
$$
I(R)=\{\gamma\in\inv{(\M_3(K)\otimes
R)}\ |\ \gamma\gamma^{\dagger}\in R\}
$$
for any (commutative) $\Q$-algebra $R$.

\vspace{1ex}
\statenumber\label{lem-unitary}{\bf .\ Lemma.}\ {\sl We have the
following isomorphisms}: 

1)\ $I(\R)\cong\mathrm{GU}(3)$. {\sl In particular,
$\ad{I}(\R)$ is a compact Lie group.}

2)\ $I(\Q_2)\cong\{(x,y)\in\GL_3(\Q_2)\times
\opp{\GL_3(\Q_2)}\,|\,xy\in\Q_2\}$. {\sl Hence in particular we have
$\ad{I}(\Q_2)\cong\PGL_3(\Q_2)$.}

\vspace{1ex}
\pf
1) is clear by definition. To prove 2) we first note the isomorphism
$$
\M_3(K)\otimes\Q_2\cong\M_3(\Q_2)\times\opp{\M_3(\Q_2)}
$$
such that the involution $\ast$ acts on the right-hand side by
interchanging the factors, i.e.\ $(x,y)^{\ast}=(y,x)$.
We have $(x,y)^{\dagger}=(HyH^{-1},H^{-1}xH)$.
Then 2) follows easily from the existence of the
decomposition (\ref{para-recall}.2).
\qed

\vspace{1ex}
\statenumber\label{para-level}{\bf .}\ 
Let $\O_K$ (resp.\ $\O_L$) be the integer ring of $K$ (resp.\ $L$).
We consider $\O_L$ as a free $\O_K$-module of rank $3$.
For a finite set $S$ of prime numbers we denote by $\A^S_{\mathrm{f}}$
the prime-to-$S$ part of the finite adele ring of $\Q$, and set 
$$
\widehat{\Z}^S=\prod_{\ell\not\in S}\Z_{\ell}.
$$
For $S$ as above and a prime $p$ we define the maximal open compact
subgroups 
\begin{center}
$\ \ C^S_{\mathrm{max}}=\{\gamma\in I(\A^S_{\mathrm{f}})\,|\,
\gamma(\O_L\otimes_{\Z}\widehat{\Z}^S)\subseteq
\O_L\otimes_{\Z}\widehat{\Z}^S\}$,\\

\vspace{1ex}
$C^{\mathrm{max}}_p=\{\gamma\in I(\Q_p)\,|\,
\gamma(\O_L\otimes_{\Z}\Z_p)\subseteq
\O_L\otimes_{\Z}\Z_p\}$
\end{center}

\noindent
of $I(\A^S_{\mathrm{f}})$ and $I(\Q_p)$, respectively.

We are going to define an open compact subgroup $C_7$ of
$C^{\mathrm{max}}_7$.
Since the prime $7$ ramifies in $K$ we have 
$$
\O_L\otimes_{\Z}\Z_7\cong\til{\Z}_7\cdot 1\oplus\til{\Z}_7\cdot\zeta
\oplus\til{\Z}_7\cdot\zeta^2,
\leqno{(\ref{para-level}.1)}
$$
where $\til{\Z}_7$ is the ramified quadratic extension of $\Z_7$.
We can therefore consider the modulo $\sqrt{-7}$ reduction of
$C^{\mathrm{max}}_7$, which takes values in $\GL_3(\F_7)$.
Now the matrix $(H$ mod $\sqrt{-7})$ is of rank $1$ and has $2$
dimensional null space, to which the action of $C^{\mathrm{max}}_7$ can
be restricted.
We therefore obtain a homomorphism 
$$
\varpi\colon C^{\mathrm{max}}_7\longrightarrow\GL_2(\F_7).
\leqno{(\ref{para-level}.2)}
$$
We take a $2$-Sylow subgroup $P$ of the subgroup 
$\GL_2(\F_7)\bigcap\{\gamma\,|\det\gamma=\pm 1\}$ 
and put $C_7=\varpi^{-1}(P)$.

\vspace{1ex}
\statenumber\label{def-mumfordgroup}{\bf .\ Definition.}\ Set
$C^2_{\mathrm{Mum}}=C^{2,7}_{\mathrm{max}}C_7$, which is an open compact
subgroup of $I(\A^2_{\mathrm{f}})$.
We define the subgroup $\Gamma_{\mathrm{Mum}}$ of
$\ad{I}(\Q_2)\cong\PGL_3(\Q_2)$ to be the image of 
\begin{center}
$I(\Q)\bigcap\,(I(\Q_2)\times C^2_{\mathrm{Mum}})
\longrightarrow\ad{I}(\Q_2)$
\end{center}
induced by the first projection, where the intersection on the left-hand 
side is taken in $I(\A_{\mathrm{f}})$.

\vspace{1ex}
\statenumber\label{thm-mumford}{\bf .\ Theorem (Mumford \cite{Mu79}).}\
{\sl The subgroup $\Gamma_{\mathrm{Mum}}$ is a uniform lattice in
$\PGL_3(\Q_2)$ which acts transitively on the vertices of the
Bruhat-Tits building attached to $\PGL_3(\Q_2)$. Moreover, the quotient 
$$
{\cal X}_{\mathrm{Mum}}=
\Gamma_{\mathrm{Mum}}\backslash\Omega^3_{\Q_2}
\leqno{(\ref{thm-mumford}.1)}
$$
of the Drinfeld symmetric space $\Omega^3_{\Q_2}$ is algebraized to a
fake projective plane $X_{\mathrm{Mum}}$.} 

\vspace{1ex}
\pf
We need to show that the subgroup $\Gamma_{\mathrm{Mum}}$ coincides with
the one constructed by Mumford in \cite{Mu79}. 
The Mumford's group is the image of 
\begin{center}
$I_1(\Q)\bigcap\,(I(\Q_2)\times C^2_{\mathrm{Mum}})
\longrightarrow\ad{I}(\Q_2)$,
\end{center}
where $I_1(\Q)=\{\gamma\in I(\Q)\,|\,\gamma\gamma^{\dagger}=1\}$.
Hence it suffice to see that every element $\gamma\in
I(\Q)\bigcap\,(I(\Q_2)\times C^2_{\mathrm{Mum}})$ 
can be written as $\gamma=k\theta$ with $k\in\inv{\Q}$
and $\theta\in I_1(\Q)\bigcap\,(I(\Q_2)\times C^2_{\mathrm{Mum}})$.
Set 
$$
k=(\gamma\gamma^{\dagger})^{-1}\det\gamma
$$
and $\theta=k^{-1}\gamma$.
Since $(\gamma$ mod $\sqrt{-7})\in P$ and $(\dagger$ mod $\sqrt{-7})$ is
trivial we have $(k$ mod $\sqrt{-7})\in P$.
Hence we have $(\theta$ mod $\sqrt{-7})\in P$.
The property $\theta\theta^{\dagger}=1$ follows from the fact that 
for any $\gamma\in I(\Q)$ we have
$\det\gamma\ovl{\det\gamma}=(\gamma\gamma^{\dagger})^3$.
\qed

\vspace{2ex}
\sectionnumber\label{section-relation}\ 
{\bf Relation between {\boldmath $X_{\mathrm{Mum}}$} and $\mathbf{Sh}$}

\vspace{1ex}
\statenumber\label{pro-inner}{\bf .\ Proposition.}\ 
{\sl Let $I$ and $G$ be the $\Q$-algebraic groups defined in
{\rm \ref{def-twistgroup}} and {\rm \ref{def-group}}, respectively.
Then we have the following$:$ 

{\rm (a)}\ \,$I$ is an inner form of $G$.

{\rm (b)}\ $I(\A^2_{\mathrm{f}})\cong G(\A^2_{\mathrm{f}})$.

\noindent
Moreover, the algebraic group $I$ is uniquely determined up to
$\Q$-isomorphisms by {\rm (a)}, {\rm (b)} and the following condition$:$

{\rm (c)}\ \,$I(\Q_2)\cong\{(x,y)\in\GL_3(\Q_2)\times
\opp{\GL_3(\Q_2)}\,|\,xy\in\Q_2\}$. }

\vspace{1ex}
\pf
First we note that, since $G(\Q_2)\not\cong I(\Q_2)$, the condition (a),
(b) and (c) imply that $\ad{I}(\R)$ is a compact Lie group.
By this the uniqueness of $I$ will follow if we prove that the group $G$
satisfies the Hasse principle; for this 
we refer the argument by Kottwitz \cite[\S 7]{Ko92}.
Our situation belongs to the case (A) in the notation as in \cite[\S
5]{Ko92}.
Kottwitz showed that validity of the Hasse principle for $G$
depends only on that for the center $Z(G)$ of $G$.
In our case we have $Z(G)=\mathrm{Res}_{K/\Q}\G_{\mathrm{m}}$, which
clearly satisfies the Hasse principle, since $K$ is a quadratic
extension of $\Q$.

Next we prove that $I$ is an inner form of $G$.
It suffices to show that the $\Q$-algebra with involution
$(\M_3(K),\dagger)$ is an inner form of $(D,\bigstar)$.
Set 
$$
R=\M_3(K)\otimes L,
$$
and consider the action $\rho$ of the Galois group $\mathrm{Gal}(L/\Q)$
on $R$ via the second factor.
The invariant part $R^{\rho}$ is isomorphic to $\M_3(K)$.
Note that the Galois group $\mathrm{Gal}(L/\Q)$ is generated by $\sigma$
and complex conjugation (cf.\ \ref{para-recall}). 
Now we consider the matrix
$$
Q=\Phi(\Pi)=\left[
\begin{array}{ccc}
&&\mu\\
1\\
&1
\end{array}
\right]\in\M_3(K),
$$
and define another Galois action $\rho'$ by $\rho'(\textrm{complex
conj.})=\rho(\textrm{complex conj.})$ and $\rho'(\sigma)(A\otimes z)=
(Q^{-1}AQ)\otimes z^{\sigma}$ for $A\in\M_3(K)$ and $z\in L$.

\vspace{1ex}
{\bf \ref{pro-inner}.1.}\ {\bf Claim.}\ {\sl The invariant part
$R^{\rho'}$ is isomorphic to $D$.}

\vspace{1ex}
Looking at the isomorphism $R\cong\M_3(L)\times\opp{\M_3(L)}$ such that
$\sigma$ (resp.\ complex conjugation) acts on $\M_3(L)\times\opp{\M_3(L)}$
diagonally (resp.\ by interchanging the factors), we deduce that
$R^{\rho'}$ is isomorphic to 
$$
D'=\{A\in\M_3(L)\,|\,Q^{-1}A^{\sigma}Q=A\}.
$$
Obviously $D'$ contains $\Phi(D)$. Comparing $\Q$-ranks we have 
$D'=\Phi(D)$, thereby the claim.

\vspace{1ex}
{\bf \ref{pro-inner}.2.}\ {\bf Lemma.}\ {\sl Let $A$ be a ring.
For $u\in\inv{A}$ with $u^{\ast}u^{-1}\in Z(A)$ we define the involution
$i_u$ on the ring $A\times\opp{A}$ by 
$$
(x,y)^{i_u}=(uyu^{-1},u^{-1}xu)
$$
for $x,y\in A$.
Then the ring with involution $(A\times\opp{A},i_u)$ is isomorphic to 
$(A\times\opp{A},i_1)$.}

\vspace{1ex}
This is well-known. In fact, the conjugation by $(1,u)$ gives the desired 
isomorphism. 

By the lemma we deduce in particular that the $L$-algebras
with involution $(R,i_H)$ and $(R,i_{\Phi(b)})$ are isomorphic to each
other.
Since $H\in\M_3(K)$, the involution $i_H$ commutes with the Galois action 
$\rho$, and hence it induces an involution on $R^{\rho}=\M_3(K)$, which
coincides with $\dagger$.
Also, since $\Phi(b)\in\M_3(K)$ and $\Phi(b)Q=Q\Phi(b)$, the involution
$i_{\Phi(b)}$ commutes with the Galois action $\rho'$ and can be
restricted to $R^{\rho'}=D$, which is nothing but $\bigstar$.
We therefore deduce that $(\M_3(K),\dagger)$ is an inner form of
$(D,\bigstar)$, thereby (a).

Finally we prove (b).
Let $\ell\neq 2$ be a rational prime.
By \ref{para-algebra}.3 we know that
$D_{\ell}\cong\M_3(K)\otimes\Q_{\ell}$.
Hence it suffices to show that the involutions $\dagger$ and $\bigstar$
are $\Q_{\ell}$-isomorphic to each other.
If the prime $\ell$ splits in $K$, the assertion follows from
\ref{pro-inner}.2.
Suppose that $\ell$ does not split in $K$.
In order to show that $\dagger$ and $\bigstar$ are isomorphic it
suffices to show that the hermitian forms $H$ and
$H'=(\lambda-\ovl{\lambda})\Phi(b)$ are isomorphic.
We can do this by the well-known fact (cf.\ \cite[App.\ 2]{MH73}): 
An isomorphism classe of hermitian forms on a fixed vector space over a
local field $F$ $($with involution $\ast)$ is determined by the residue
class of determinant modulo $\mathrm{N}_{F/F_0}(\inv{F})$, where $F_0$
is the fixed field of $\ast$.
We know that
$\det H=7$ and $\det H'=7^2$.
If $\ell\neq 7$ then $K_{\ell}/\Q_{\ell}$ is an unramified quadratic
extension. 
Since both $\det H$ and $\det H'$ are unit integer and since every unit
integer of $\Q_{\ell}$ is a norm of an element in the unramified quadratic 
extension, we have the desired result in this case.
If $\ell=7$ then $K_7=\Q_7(\sqrt{-7})$, and hence both $\det H$ and
$\det H'$ are obviously norms, thereby the result.

Now we have proved (b), and hence the proof of the proposition is
finished.
\qed

\vspace{1ex}
\statenumber\label{para-uniformization}{\bf .}\ 
We fix an isomorphism $\varphi\colon I(\A^2_{\mathrm{f}})\cong
G(\A^2_{\mathrm{f}})$ 
as in \ref{pro-inner} (b).
Let $C^{\mathrm{max}}_2\subset G(\Q_2)$ be a maximal open compact
subgroup defined by 
$$
C^{\mathrm{max}}_2=\{\gamma\in G(\Q_2)\,|\,
(\O_D\otimes_{\Z}\Z_2)\gamma\subseteq
\O_D\otimes_{\Z}\Z_2\}.
$$
We define an open compact subgroup $C\subset G(\A_{\mathrm{f}})$ by 
$$
C=C^{\mathrm{max}}_2\varphi(C^2_{\mathrm{Mum}})
\leqno{(\ref{para-uniformization}.1)}
$$
(cf.\ \ref{def-mumfordgroup}).
Let us fix an embedding $\nu\colon\ovl{\Q}\hookrightarrow\ovl{\Q}_2$.
Let $E_{\nu}$ be the $\nu$-adic completion of the Shimura field $E$
$(=\varepsilon(K))$, and $\breve{E}_{\nu}$ the maximal unramified
extension of $E_{\nu}$.
By \cite[6.51]{RZ96} (see also
\cite[1.12]{BZ95}\cite{Va98a}\cite{Va98b}) and the uniqueness of $I$ as
in \ref{pro-inner} we have an isomorphism of formal schemes
$$
I(\Q)\backslash(\widehat{\Omega}^2_{\Q_2}\times_{\Spf\Z_2}
\Spf\O_{\breve{E}_{\nu}})\times 
G(\A_{\mathrm{f}})/C\stackrel{\sim}{\longrightarrow}
\mathrm{Sh}^{\wedge}_C\times_{\Spf\O_{E_{\nu}}}
\Spf\O_{\breve{E}_{\nu}},
\leqno{(\ref{para-uniformization}.2)}
$$
where $\mathrm{Sh}^{\wedge}_C$ is the formal completion of a model of 
$\mathrm{Sh}_C$ over $\Spec\O_{E_{\nu}}$ along the closed fiber. 

\vspace{1ex}
\statenumber\label{thm-main}{\bf .\ Theorem.}\ {\sl The canonical model
$\mathrm{Sh}_C$ of the Shimura variety ${\cal S}h_C$ has exactly three
geometrically connected components defined over $\Q(\zeta_7)$ which are
permuted by the action of $\mathrm{Gal}(\Q(\zeta_7)/\Q(\sqrt{-7}))$.
Moreover, for any connected component $S$ of $\mathrm{Sh}_C$ there
exists an isomorphism of schemes
$$
X_{\mathrm{Mum}}\times_{\Spec\Q_2}\Spec\Q_2(\zeta_7)
\stackrel{\sim}{\longrightarrow} 
S\times_{\Spec\Q(\zeta_7)}\Spec\Q_2(\zeta_7)
\leqno{(\ref{thm-main}.1)}
$$
over $\Q_2(\zeta_7)$, where $X_{\mathrm{Mum}}$ is the Mumford's fake
projective plane $($cf.\ {\rm \ref{thm-mumford}}$)$.}

\vspace{1ex}
\pf
First we note that by \cite[6.51]{RZ96} we can take the isomorphism as
in (\ref{para-uniformization}.2) such that the natural descent datum on
the right-hand side induces on the left-hand side the natural descent
datum on the first factor multiplied with the action of 
$$
(\Pi,\Pi^{-1})\in
G(\Q_2)\subset\inv{D}_{\lambda}\times\opp{\inv{D}_{\lambda}} 
$$
on $G(\A_{\mathrm{f}})/C$.
Since $\Pi\in\O_{D_{\lambda}}$ and
$\Pi^{-1}=\ovl{\mu}\Pi^2\in\O_{D_{\ovl{\lambda}}}$ we have
$(\Pi,\Pi^{-1})\in C^{\mathrm{max}}_2$.
Therefore the action of $(\Pi,\Pi^{-1})$ is trivial on
$G(\A_{\mathrm{f}})/C$.
Hence the isomorphism (\ref{para-uniformization}.2) descends to an
isomorphism 
$$
I(\Q)\backslash\widehat{\Omega}^2_{\Q_2}\times 
G(\A_{\mathrm{f}})/C\stackrel{\sim}{\longrightarrow}
\mathrm{Sh}^{\wedge}_C
\leqno{(\ref{thm-main}.2)}
$$
over $\Z_2\cong\O_{E_{\nu}}$. 

Next we note that the left-hand side of (\ref{para-uniformization}.2) is
the disjoint union of formal schemes of form 
$$
\Gamma_{\gamma}\backslash(\widehat{\Omega}^2_{\Q_2}\times_{\Spf\Z_2}
\Spf\O_{\breve{E}_{\nu}}),
$$
where $\Gamma_{\gamma}=$ the image of $I(\Q)\bigcap(I(\Q_2)\times(\gamma
C^2_{\mathrm{Mum}}\gamma^{-1}))$ in $\ad{I}(\Q)$ for
$\gamma\in\G(\A^2_{\mathrm{f}})$. 
Note also that $\Gamma_1=\Gamma_{\mathrm{Mum}}$ (cf.\
\ref{def-mumfordgroup}).

By these observations it suffices to prove that the canonical model
$\mathrm{Sh}_C$ of the Shimura variety has three connected components 
as in the theorem.
To do this we consider the torus $T$ over $\Q$ defined by
$$
T=\{(k,f)\in\mathrm{Res}_{K/\Q}\G_{\mathrm{m},K}
\times\G_{\mathrm{m},\Q}\ |\ k\ovl{k}=f^3\}\cong\mathrm{Res}_{K/\Q}
\G_{\mathrm{m},K}
$$
and define an epimorphism $\vartheta\colon G\rightarrow T$ by 
$\gamma\mapsto(\mathrm{N}^{\circ}_{\opp{D}/K}\,\gamma,c(\gamma))$, where
$\mathrm{N}^{\circ}_{\opp{D}/K}$ is the reduced norm (the last $\cong$
is given by $(k,f)\mapsto kf^{-1}$).
The kernel of $\vartheta$ is the derived group of $G$.
By \cite[\S 2]{De71} we know that the morphism $\vartheta$ induces the
canonical morphism of Shimura varieties
$$
{\cal S}h_C=G(\Q)\setminus X_{\infty}\times 
G(\A_{\mathrm{f}})/C
\longrightarrow T(\Q)_{+}\setminus T(\A_{\mathrm{f}})/
\vartheta(C)
$$
which induces the Galois-equivariant bijection between the sets of all 
connected components of canonical models of each side.

To calculate $\vartheta(C)$ we have only to look at the component
$\varphi(C_7)$ at $7$.

\vspace{1ex}
{\bf \ref{thm-main}.3.}\ {\bf Claim.}\ {\sl Under the canonical
identification $T(\Q_7)\cong\til{\Q}^{\times}_7$, where
$\til{\Q}^{\times}_7$ is the ramified quadratic extension of $\Q_7$, we
have} 
$$
\vartheta(\varphi(C_7))\cong\pm 1+\sqrt{-7}\,\til{\Z}^{\times}_7,
$$
{\sl where $\til{\Z}^{\times}_7$ is the integer ring of
$\til{\Q}^{\times}_7$.} 

\vspace{1ex}
First we note that the map $\vartheta\circ\varphi$ is described by
$$
\vartheta\circ\varphi\colon
C_7\ni\gamma\mapsto(\gamma\gamma^{\dagger})^{-1}\det\gamma
\in\til{\Z}^{\times}_7.
$$
Since $\pm 1+7\,\inv{\Z}_7\subset C_7$ we immediately deduce that 
$\vartheta(\varphi(C_7))$ contains $\pm 1+7\,\inv{\Z}_7$.
Moreover, since $\pm 1+\sqrt{-7}\,\til{\Z}^{\times}_7\subset C_7$,
$\vartheta(\varphi(C_7))$ contains the element $k^2\ovl{k}^{-1}$ for any 
$k\in\pm 1+\sqrt{-7}\,\til{\Z}^{\times}_7\subset C_7$; but, since 
$k\ovl{k}\in\vartheta(\varphi(C_7))$, we deduce that 
$$
(\pm 1+\sqrt{-7}\,\til{\Z}^{\times}_7)^3=\pm
1+\sqrt{-7}\,\til{\Z}^{\times}_7 
\subseteq
\vartheta(\varphi(C_7)).
$$
Hence the claim follows if we check that $(\vartheta(\varphi(C_7))$ mod
$\sqrt{-7})=\{\pm 1\}$.
By an argument similar to that in the proof of \ref{thm-mumford} any
$\gamma\in C_7$ has a decomposition $\gamma=k\theta$ in $C_7$ with
$\theta\theta^{\dagger}=1$.
Since $\det\theta\det\theta^{\dagger}\equiv(\det\theta)^2\equiv 1$
mod $\sqrt{-7}$ we have $\vartheta(\varphi(\gamma))\equiv\pm k$.
But, since $k=\gamma\theta^{-1}$ mod $\sqrt{-7}$ belongs to $P$, the
order of $k$ is either one or two; hence $k\equiv\pm 1$ mod $\sqrt{-7}$.
We therefore get $\vartheta(\varphi(\gamma))\equiv\pm 1$ mod
$\sqrt{-7}$ for any $\gamma\in C_7$, which proves the claim.

Now in view of the theory of complex multiplication we find that the
Shimura variety 
$T(\Q)_{+}\backslash T(\A_{\mathrm{f}})/\vartheta(C)$ has three
connected components having field of definition $L=\Q(\zeta)$ which are
permuted by $\mathrm{Gal}(L/K)$, since the ray class field corresponding 
to $C$ is easily checked to be $L$.
We therefore proved the all statements in the theorem.
\qed

\vspace{1ex}
\statenumber\label{cor-main}{\bf .\ Corollary.}\ {\sl 
Each connected component of the complex surface ${\cal S}h_C$ is a fake
projective plane, 
which has a field of definition $\Q(\zeta_7)$.
Moreover, it is an arithmetic quotient of complex
unit-ball.}
\qed

\vspace{1ex}
{\it Acknowledgments.}\ 
The author's viewpoint of this work owes much to suggestions 
by Professor Yves Andr\'e. 
The author expresses gratitude to him for valuable discussions.
The author thanks Professor Yoichi Miyaoka for his interest in this work 
and several valuable suggestions on the manuscript. 
The author thanks Max-Planck-Institut f\"ur Mathematik Bonn for the
nice hospitality.

\begin{small}

{\sc Graduate School of Mathematics, Kyushu University, 6-10-1\ Hakozaki
Higashi-ku, Fukuoka 812, Japan.}

e-mail: {\tt fkato@math.kyushu-u.ac.jp}
\end{small}

\begin{thebibliography}{99}
\bibitem[An98]{An98}
              Andr\'e, Y.: {\it $p$-adic orbifolds and $p$-adic triangle 
              groups}, preprint, Institut de Math.\ Jussieu (1998).
\bibitem[Ba99]{Ba99}
              Barlow, R.N.: {\it Zero-cycles on Mumford's surface},
              Math.\ Proc.\ Cambridge Philos.\ Soc.\ {\bf 126} (1999),
              no.\ 3, 505--510.  
\bibitem[BZ95]{BZ95}
              Boutot, J-F., Zink, Th.: {\it The $p$-adic Uniformization
              of Shimura Curves}, preprint 95--107, Universit\"at
              Bielefeld (1995).
\bibitem[\v{C}e76]{Ce76}
              \v{C}erednik, I.V.: {\it Uniformization of algebraic
              curves by discrete arithmetic subgroups of 
              $\mathrm{PGL}_2(k_{w})$ with compact quotient spaces}, 
              Mat.\ Sd., {\bf 100} (1976), 59--88.
\bibitem[De71]{De71}
              Deligne, P.: {\it Travaux de Shimura}, S\'em.\ Bourbaki
              Expos\'e 389, Lect.\ Notes in Math. {\bf 244} (1971),
              Springer, Heidelberg.
\bibitem[Dr76]{Dr76}
              Drinfeld, V.G.: {\it Coverings of $p$-adic symmetric
              regions}, Funct.\ Anal.\ and Appl.\ {\bf 10} (1976), 
              29--40.
\bibitem[IK98]{IK98}
              Ishida, M., Kato, F.: {\it The strong rigidity theorem for
              non-archimedean uniformization}, T\^ohoku Math.\ J., 
              {\bf 50} (1998), 537--555.
\bibitem[Ko92]{Ko92}   
              Kottwitz, R.E.: {\it Points on some Shimura varieties over
              finite fields}, Journal Amer.\ Math.\ Soc.\ {\bf 5}
              (1992), 373--444.
\bibitem[Ku79]{Ku79}
              Kurihara, A.: {\it On some examples of equations defining
              Shimura curves and the Mumford uniformization}, J.\ Fac.\ 
              Sci.\ Uni.\ Tokyo, Sec.\ IA {\bf 25} (1979), 277--300.
\bibitem[Ku94]{Ku94}
              Kurihara, A.: {\it On $p$-adic Poincar\'e series and
              Shimura curves}, Internat.\ J.\ Math.\ {\bf 5} (1994),
              747--763. 
\bibitem[MH73]{MH73}
              Milnor, J., Husemoller, D.: {\it Symmetric Bilinear
              Forms}, Ergebnisse der Mathematik und ihrer Grenzgebiete, 
              Band 73, Springer-Verlag, Berlin, Heidelberg, New York
              (1973).
\bibitem[Mu79]{Mu79}
              Mumford, D.: {\it An algebraic surface with $K$ ample, 
              $(K^2)=9,\ \mbox{\rm p}_{\mbox{\rm g}}=\mbox{\rm q}=0$},
              in Contribution to Algebraic Geometry, 
              Johns Hopkins Univ.\ Press (1979), 233--244.
\bibitem[RZ96]{RZ96}
              Rapoport, M., Zink, Th.: {\it Period Spaces for
              $p$-divisible Groups}, Annals of Math.\ Studies {\bf 141}
              Princeton University Press (1997).
\bibitem[Va98a]{Va98a}
              Varshavsky, Y.: {\it $p$-adic uniformation of unitary
              Shimura varieties}, Inst.\ Hautes\ \'Etudes Sci.\ Publ.\
              Math.\ {\bf 87} (1998), 57--119. 
\bibitem[Va98b]{Va98b}
              Varshavsky, Y.: {\it $p$-adic uniformization of unitary
              Shimura varieties.\ II}, J.\ Differential Geom.\ {\bf 49}
              (1998), no.\ 1, 75--113. 
\end{thebibliography}
\end{document}